\documentclass[12pt,reqno]{amsart}
\usepackage[english,russian]{babel}
\usepackage{amssymb,amsmath,amsthm}
\usepackage{graphicx}
\usepackage{epstopdf}
\usepackage{cite}
\usepackage{setspace}
\usepackage{tikz}
\newcommand{\bc}{\begin{center}}
\newcommand{\ec}{\end{center}}
\newcommand{\be}{\begin{equation}}
\newcommand{\ee}{\end{equation}}
\newcommand{\ba}{\begin{array}}
\newcommand{\ea}{\end{array}}
\newcommand{\edc}{\end{document}}

\begin{document}

\title[Hard-Core model]
{New class of Gibbs measures for two state Hard-Core model on a Cayley tree}

\author{R.~M.~Khakimov, M~T.~Makhammadaliev and F.~H.~Haydarov}

\address{R.~M.~Khakimov, M.~T.~Makhammadaliev \\ V. I. Romanovskiy Institute of Mathematics of the Academy of Sciences of Uzbekistan, Tashkent, Uzbekistan. \\ Namangan State University, Namangan, Uzbekistan.}
 \email {rustam-7102@rambler.ru, mmtmuxtor93@mail.ru}

\address{F.~H.~Haydarov \\ V. I. Romanovskiy Institute of Mathematics of the Academy of Sciences of Uzbekistan, Tashkent, Uzbekistan.}
\address{AKFA University, Tashkent, Uzbekistan.}
\address{National University of Uzbekistan, Tashkent , Uzbekistan.}
\email {haydarov\_imc@mail.ru}

 \selectlanguage{english} \begin{abstract} In this paper, we consider a Hard-Core $(HC)$ model with two spin values on Cayley trees. The conception of alternative Gibbs measure is introduced and translational invariance conditions for alternative Gibbs measures are found.  Also, we show that the existence of alternative Gibbs measures which are not translation-invariant. In addition, we study free energy of the model.
\end{abstract}

\maketitle

{\bf{Key words.}}  Cayley tree, configuration, hard-core model, Gibbs measure, translation-invariant measure, Alternating Gibbs Measure, free energy.\\

AMS Subject Classification: 20B07, 20E06.\\

\section{Introduction}
The problems arising in the study of the thermodynamic properties of physical and biological systems are typically solved within the framework of the theory of Gibbs measures. The Gibbs measure is a fundamental concept that determines the probability of a microscopic state of a given physical system (defined by a specific Hamiltonian). It is known that each Gibbs measure is associated with one phase of a physical system, and if the Gibbs measure is not unique, then there exists a phase transition. For a wide class of Hamiltonians, it is known that the set of all Gibbs measures (corresponding to a given Hamiltonian) is a nonempty, convex, compact subset of the set of all probability measures (see, e.g., \protect{\cite{6}, \cite{111}}) and each point of this convex set can be uniquely decomposed in terms of its extreme
points. In this regard, it is of particular interest to describe all the extreme points of this convex set, i.e., extreme Gibbs measures.

For convenience, we first describe the basic concepts used in this paper and then give the statement of the problem and the history of its study.

\textbf{The Cayley tree.} Let $\Im^k=(V,L,i)$, $ k\geq 1 $, be the Cayley tree of order $k$, i.e., an infinite tree with exactly $k+1$ edges coming out of each vertex, and let $V$ be the set of vertices, $L$  the set of edges of $\Im^k$ and $i$ is the incidence function setting each edge $l\in L$ into correspondence with its endpoints $x, y \in V$. If $i(l)=\{x, y\}$, then the vertices $x$ and $y$ are called the {\it nearest neighbors}, denoted by $l = \langle x, y \rangle $.

For an arbitrary point $x^0\in V$ we set
$$  W_n = \ \{x\in V\ \ | \ \ d (x,
x^0) =n \}, \ \ V_n = \bigcup_{m=0}^{n}{W_m},
\ \ L_n = \ \{l =
\langle x, y\rangle \in L \ \ | \ \ x, y \in V_n \},
$$
where $d(x,y), x, y \in V$ is the distance between $x$ and $y$ on the Cayley tree, i.e., the number of edges of the path connecting $x$ and $y$.

The set of the direct successors of $x$ is denoted by $S(x)$, i.e., if $x\in W_n$, then
$$S(x)=\{y_i\in W_{n+1} |  d(x,y_i)=1, i=1,2,\ldots, k \}.$$

\textbf{The HC-model.} We assume that $\Phi=\{0, 1 \},$  and $\sigma\in \Omega=\Phi^V$ is a configuration, i.e., $\sigma=\{\sigma(x)\in \Phi: x\in V\}$, where $\sigma (x)=1$ means that
the vertex $x$ on the Cayley tree is \textit{occupied}, and $\sigma(x)=0$
means it is \textit{vacant}. The configuration $\sigma$ is said to be an \textit{admissible} if $\sigma(x)\sigma(y)=0$ for any neighboring
$\langle x,y\rangle$ from $V$ ($V_n$ or $W_n$, respectively) and we denote the set of such configurations by $\Omega$ ($\Omega_{V_n}$ and $\Omega_{W_n})$. Obviously,
$\Omega\subset\Phi^V.$

Concatenation configurations $\sigma_{n-1}\in\Phi^{V_{n-1}}$ and $\omega_n\in\Phi^{W_{n}}$ is defined by the following formula (see \protect{\cite{1111}})
$$
\sigma_{n-1}\vee\omega_n=\{\{\sigma_{n-1}(x), x\in V_{n-1}\},
\{\omega_n(y), y\in W_n\}\}.
$$

The Hamiltonian of the \emph{HC}-model is defined by the formula
$$H(\sigma)=\left\{
\begin{array}{ll}
    J \sum\limits_{x\in{V}}{\sigma(x),} \ \ \ $ if $ \sigma \in\Omega $,$ \\
   +\infty ,\ \ \ \ \ \ \ \ \ \ $  \ if $ \sigma \ \notin \Omega $,$ \\
\end{array}
\right. $$ where $J\in R$.

\textbf{Finite-dimensional distributions.} Let $\mathbf{B}$ be the $\sigma$-algebra generated by
cylindrical sets with finite base of $\Omega$. For any $n$ we let $\mathbf{B}_{V_n}=\{\sigma\in\Omega:
\sigma|_{V_n}=\sigma_n\}$ denote the subalgebra of $\mathbf{B},$ where
$\sigma|_{V_n}$-restriction of $\sigma$ to $V_n$ and $\sigma_n: x\in V_n
\mapsto \sigma_n(x)$ an admissible configuration in $V_n.$

\textbf{Definition 1}. For $\lambda >0$, the \emph{HC}-model Gibbs measure is
a probability measure $\mu$ on $(\Omega , \textbf{B})$ such that for
any $n$ and $\sigma_n\in \Omega_{V_n}$
$$
\mu \{\sigma \in \Omega:\sigma|_{V_n}=\sigma_n\}=
\int_{\Omega}\mu(d\omega)P_n(\sigma_n|\omega_{W_{n+1}}),
$$
where
$$
P_n(\sigma_n|\omega_{W_{n+1}})=\frac{e^{-H(\sigma_n)}}{Z_{n}
(\lambda ; \omega |_{W_{n+1}})}\textbf{1}(\sigma_n \vee \omega
|_{W_{n+1}}\in\Omega_{V_{n+1}}).
$$
Here $Z_n(\lambda ; \omega|_{W_{n+1}})$- normalization factor with
boundary condition $\omega|_{W_n}$:
$$
Z_n (\lambda ; \omega|_{W_{n+1}})=\sum_{\widetilde{\sigma}_n \in
\Omega_{V_n}}
e^{-H(\widetilde{\sigma}_n)}\textbf{1}(\widetilde{\sigma}_n\vee
\omega|_{W_{n+1}}\in \Omega_{V_{n+1}}).
$$
For $\sigma_n\in\Omega_{V_n}$ we denote that
$$\#\sigma_n=\sum\limits_{x\in V_n}{\mathbf 1}(\sigma_n(x)\geq 1)$$
which the number of occupied vertices in $\sigma_n$.

Let $z:x\mapsto z_x=(z_{0,x}, z_{1,x}) \in R^2_+$ vector-valued function on $V$.
For $n=1,2,\ldots$ and $\lambda>0$ consider the probability measure $\mu^{(n)}$ on $\Omega_{V_n}$,
defined as
\begin{equation} \label{e1}
\mu^{(n)}(\sigma_n)=\frac{1}{Z_n}\lambda^{\#\sigma_n} \prod_{x\in
W_n}z_{\sigma(x),x}.
\end{equation}
Here $Z_n$ is the normalizing divisor:
$$Z_n=\sum_{{\widetilde\sigma}_n\in\Omega_{V_n}}
\lambda^{\#{\widetilde\sigma}_n}\prod_{x\in W_n}
z_{{\widetilde\sigma}(x),x}.$$

The sequence of probability measures $\mu^{(n)}$ is said to be consistent if for any $n\geq1$ and
$\sigma_{n-1}\in\Omega_{V_{n-1}}$:

\begin{equation} \label{e2}
\sum_{\omega_n\in\Omega_{W_n}}
\mu^{(n)}(\sigma_{n-1}\vee\omega_n){\mathbf 1}(
\sigma_{n-1}\vee\omega_n\in\Omega_{V_n})=
\mu^{(n-1)}(\sigma_{n-1}).
\end{equation}

In this case, there is a unique measure $\mu$ on $(\Omega,\textbf{B})$ such that for all $n$ and $\sigma_n\in \Omega_{V_n}$
$$\mu(\{\sigma|_{V_n}=\sigma_n\})=\mu^{(n)}(\sigma_n).$$

\textbf{Definition 2.} The measure $\mu$ that is the limit of a sequence $\mu^{(n)}$ defined by (\ref{e1}) with consistency condition (\ref{e2}) is called the splitting \emph{HC}-Gibbs measure (SGM) with $\lambda>0$ corresponding to the function $z:\,x\in V\setminus\{x^0\}\mapsto z_x$.
Moreover, an \emph{HC}-Gibbs measure corresponding to a constant function $z_x\equiv z$
is said to be translation-invariant (TI).

\textbf{Problem statement.} The main task is to study the structure of the set $\mathcal G(H)$ of all Gibbs measures corresponding to a given Hamiltonian $H$.

A measure $ \mu \in \mathcal G(H) $ is called extreme if it cannot be expressed as $ \mu = \lambda \mu_1 + (1-\lambda)\mu_2$ for some $ \mu_1, \mu_2 \in \mathcal G(H) $ with $ \mu_1 \ne \mu_2 $.

As noted above, the set $ \mathcal G(H) $ of all Gibbs measures (for a given Hamiltonian $H$) is a nonempty convex compact set $ \mathcal G (H) $ in the space of all probability measures on $ \Omega$.

Using theorem (12.6) in \protect{\cite{6}} and section 1.2.4 in \protect{\cite{BR}}, we can note the following.

\begin {itemize}
\item \ \ \emph{Any extreme Gibbs measure
$ \mu \in \mathcal G(H) $ is an SGM; therefore, the problem of describing Gibbs measures
reduces to describing the set of SGMs. For each fixed temperature, the description of the set $ \mathcal G(H) $ is equivalent to a complete description of the set of all extreme SGMs, and hence we are only interested in SGMs on the Cayley tree.}

\item \ \ \emph{Any SGM corresponds to the solution of Eq. (\ref{e3}) (see below). Thus, our main
task reduces to solving functional equation (\ref{e3})}.
\end{itemize}

It is known \protect{\cite{7}} that each Gibbs measure for \emph{HC}-model on the Cayley tree can be associated with the collection of values
$z=\{z_x,x\in V\}$ satisfying
\begin{equation} \label{e3}
z_x=\prod_{y \in S(x)}(1+\lambda z_y)^{-1},
\end{equation}
where $\lambda=e^{-J\beta}>0$ is a parameter, $\beta={1\over T}$, $T>0$ is a temperature.

Let $G_k$ be a free product of $k+1$ cyclic groups $\{e,a_i\}$ of order two with the respective generators $a_1,a_2,...,a_{k+1}, a_i^2=e$.
There is a one-to-one correspondence between the set of vertices $V$ of the Cayley
tree of order $k$ and the group $G_k$ (see \protect{\cite{11, GR2003, 1}}).

Let $\widehat{G}_k$ be a normal divisor of a finite index $r\geq 1$ and $G_k/\widehat{G}_k=\{H_1,...,H_r\}$ be the quotient group.

\textbf{Definition 3.}  A collection of quantities $z=\{z_x, x\in G_k\}$ is said to be $\widehat{G}_k$-periodic if $z_{yx}=z_x$ for $\forall x\in G_k, y\in\widehat{G}_k.$ The $G_k$-periodic collections are called translation invariant.

For any $x\in G_k $, the set $\{y\in G_k: \langle
x,y\rangle\}\setminus S(x)$ contains a unique element denoted by $x_{\downarrow}$ (see \protect{\cite{8, 5}}).

\textbf{Definition 4.}  A collection of quantities $z=\{z_x,x\in G_k\}$ is called $\widehat{G}_k$-weakly periodic if $z_x=z_{ij}$ for any $x\in H_i$, $x_\downarrow \in H_j$ for any $x \in G_k$.

\textbf{Definition 5}. A measure $\mu$ is called $\widehat{G}_k$-(weakly) periodic if it  corresponds to a $\widehat{G}_k$-(weakly) periodic collection of quantities $z$.

\textbf{History of the study of SGMs for the \emph{HC}-model.} We present a brief overview of the work related to the Potts model on the Cayley tree.

In \protect{\cite{Maz}} A. Mazel and Yu. Suhov introduced and studied the \emph{HC}-model on the $d$-dimensional lattice $ \mathbb Z^d$. Studying Gibbs measures for the two state \emph{HC}-model on the Cayley tree was the topic in \protect{\cite{7}-\cite{RKhM}}. In \protect{\cite{7}}, the uniqueness of the translation-invariant Gibbs measure and the nonuniqueness of periodic Gibbs
measures for the \emph{HC}-model were proved. For the parameters of the \emph{HC}-model, a sufficient condition was also found in \protect{\cite{7}} under which the translation-invariant Gibbs measure is nonextreme. In the case where the translation-invariant Gibbs measure is extreme, a sufficient condition was found in \protect{\cite{Mar}}. The range of the extremes of this measure was extended in \protect{\cite{RKh2}}. Weakly periodic Gibbs measures for the \emph{HC}-model in the case of a normal divisor of index 2 were studied in \protect{\cite{Kh}} and a complete description of the weakly periodic Gibbs measures was given.

Weakly periodic Gibbs measures for the \emph{HC}-model in the case of a normal divisor of index 4 were studied in \protect{\cite{Kh1}-\cite{KhRM}}. In this case conditions for the existence of weakly periodic (nonperiodic) Gibbs measures are found. We also found conditions for the translation-invariance of the weakly periodic Gibbs measures (see Chap. 7 in \protect{\cite{Rb}} for other HC model properties and their generalizations on a Cayley tree).

In this paper, we study a two-state \emph{HC}-model on a Cayley tree. The concept of an alternative Gibbs measure is introduced. Translational invariance conditions for alternative Gibbs measures are found. In addition, the existence of alternative Gibbs measures that are not translation invariant is proved.

\section{A new class of Gibbs measures}

We consider the half-tree. Namely the root $x^0$ has $k$ nearest neighbors. We construct below new solutions of the functional equation (\ref{e3}). Consider the
following matrix
$$M=\begin{pmatrix} m & k-m \\ r & k-r \end{pmatrix}$$
where  $0\leq m\leq k$ and $0\leq r\leq k$ are non-negative integers.
This matrix defines the number of times the values $h$ and $l$ occur in the set $S(x)$ for each $z_x\in \{h,l\}$. More precisely, the boundary condition
$z=\{z_x,x\in G_k\}$ with fields taking  values $h$, $l$ defined by the following steps:

$\bullet$ if at vertex $x$ we have $z_x=h$, then the function $z_y$, which gives real values to each vertex $y\in S(x)$ by the following rule
$$\begin{cases}
h ~\mbox{on} ~m~ \mbox{vertices of} ~ S(x),\\
l ~\mbox{on} ~k-m ~\mbox{remaining vertices,}
\end{cases}$$

$\bullet$ if at vertex $x$ we have $z_x=l$, then the function $z_y$, which gives real values to each vertex $y\in S(x)$ by the following rule
$$\begin{cases}
l ~\mbox{on} ~r~ \mbox{vertices of} ~ S(x),\\
h ~\mbox{on} ~k-r ~\mbox{remaining vertices.}
\end{cases}$$
For an example of such a function see Fig.1.

Then the system (\ref{e3}) has the form
\begin{equation} \label{e4}
\left\{
\begin{array}{ll}
    h=\frac{1}{\left(1+\lambda h\right)^m}\cdot \frac{1}{\left(1+\lambda l\right)^{k-m}}, \\
    l=\frac{1}{\left(1+\lambda l\right)^r}\cdot \frac{1}{\left(1+\lambda h\right)^{k-r}},  \\
    \end{array}
\right.
\end{equation}
where $l>0, h>0, \lambda>0.$

\begin{figure}
\scalebox{1}{
\begin{tikzpicture}[level distance=2.5cm,
level 1/.style={sibling distance=2.3cm},
level 2/.style={sibling distance=2.5cm},
level 3/.style={sibling distance=.3cm}]
\node {$h$} [grow'=up]
    child[sibling distance=3.5cm] {node {$l$}
        child[sibling distance=.5cm] foreach \name in {h,h,h,h,l} { node
        {$\name$} }
        }
        child[sibling distance=3.4cm] {node {$l$}
            child[sibling distance=.5cm] foreach \name in {h,h,h,h,l} { node
        {$\name$} }
        }
        child[sibling distance=2.5cm] {node {$h$}
        child[sibling distance=.5cm] foreach \name in {l,l,h,h,h} { node
        {$\name$} }
    }
    child[sibling distance=3.4cm]{node {$h$}
        child[sibling distance=.5cm] foreach \name in {l,l,h,h,h} { node
        {$\name$} }
    }
    child[sibling distance=3.5cm] {node {$h$}
        child[sibling distance=.5cm] foreach \name in {l,l,h,h,h} { node
        {$\name$} }
    }
;
\end{tikzpicture}
}
\begin{center}{\footnotesize \noindent Figure 1. In this figure the values of function $z_x$ on the vertices of the Cayley tree
of order 5 are shown. This is the case when $m=3$ and $r=2$.}
\end{center}
\end{figure}

As was mentioned above, for any boundary condition satisfying the functional equation (\ref{e3}) there exists a unique Gibbs measure. A measure constructed in this way and which is not translation-invariant is called alternative Gibbs
measure (AGM) and denoted as $\mu_{m,r}$.

\textbf{Remark 1.} Note that the solution $l=h$ in (\ref{e4}) corresponds to the only TIGM for the \emph{HC}-model (see \protect{\cite{7}}). Therefore, we are interested
in solutions of the form $l\neq h.$

\textbf{Remark 2.} From (\ref{e4}) for $m=r=0$ we obtain a system of equations whose solutions correspond to the $G^{(2)}_k$-periodic Gibbs measures for the \emph{HC}-model.

The following theorem holds.

\textbf{Theorem 1.} Let $k\geq2$. If $m+r\geq k-1$ then for the \emph{HC}-model there is a unique AGM, which coincides with the unique TIGM.

\textbf{Proof.} For convenience, we denote $h=x$ and $l=y$. Then (\ref{e4}) can be rewritten as follows:
\begin{equation} \label{e5}
\left\{
\begin{array}{ll}
    x=\frac{1}{(1+\lambda x)^m}\cdot \frac{1}{(1+\lambda y)^{k-m}}, \\
    y=\frac{1}{(1+\lambda y)^r}\cdot \frac{1}{(1+\lambda x)^{k-r}}.  \\
\end{array}
\right.
\end{equation}

If the first equation (\ref{e5}) is divided by the second, then
\begin{equation}\label{e6}
\frac xy=\left(\frac{1+\lambda x}{1+ \lambda y}\right)^{k-m-r}
\end{equation}

We denote $m+r-k=t$, $t\geq-1$. Then by (\ref{e6}) we have
$$x\left(1+ \lambda x\right)^t=y\left(1+ \lambda y\right)^t$$
It is easy to check that the function $f(x)=x\left(1+ \lambda x\right)^t$ is increasing for $t\geq-1$.
Therefore, if $m+r\geq k-1$, then the system of equations (\ref{e5}) has only a solution of the form $x=y$, and this solution corresponds to the TIGM and is known to
be unique. The theorem is proved.

By theorem 1 follows the next

\textbf{Consequence.} Let $k\geq2$. If  there are AGMs (non TI) for the \emph{HC}-model, then $m+r\leq k-2$.

Let $k-m-r=n$ $(n\in N, n\geq2)$. Then from the (\ref{e6}) we get
$$x\left(1+ \lambda y\right)^n=y\left(1+ \lambda x\right)^n.$$
From this equation after simple algebra, we obtain the equation
$$(y-x)\Big(-1+C_n^2\lambda^2xy +C_n^3\lambda^3xy (x+y)+\ldots + C_n^n\lambda^n xy(x^{n-2}+x^{n-3}y+x^{n-4}y^2+\ldots +y^{n-2})\Big)=0.$$
Hence $x=y$ or $g(x,y)=0$, where
$$
g(x,y)=C_n^2\lambda^2xy +C_n^3\lambda^3xy (x+y)+\ldots + C_n^n\lambda^n xy(x^{n-2}+x^{n-3}y+x^{n-4}y^2+\ldots +y^{n-2})-1.
$$
In the case $x=y$ the corresponding measure is TIGM.

The case $x\neq y$. We consider the equation $g(x,y)=0$ with respect to the variable $x$ (or $y$). Then it's clear that $g(0,y)=-1$ and $g(x,y)\rightarrow+\infty$ for $x\rightarrow+\infty$. Then the equation $g(x,y)=0$ for variable $x$ has at least one positive root.
On the other hand, by Descartes' theorem, the equation $g(x,y)=0$ for variable $x$ has at most one positive root. Hence the equation $g(x,y)=0$ for variable $x$ has exactly one positive root, i.e., there exists a solution $(x,y)$ of the system of equations (\ref{e5}), different from $(x,x) $. Thus, the following statement is true.

\textbf{Statement 1.} Let $k\geq2$. If $m+r\leq k-2$ then for the \emph{HC}-model there exists AGM (not TI).

In particular, if $\lambda^2xy=1$ for $m+r=k-2$ $(n=2)$ or if $\lambda^2xy(3+\lambda(x+y))=1$ for $m+r=k-3$ $(n=3)$, then in both cases there exists AGM (not TI).

The case $x=y$. We check the multiplicity of the root $x=y$. In this case from $g(x,x)=0$, we have
\begin{equation}\label{e7a}
C_n^2\lambda^2x^2 +2C_n^3\lambda^3x^3+\ldots + (n-1)C_n^n\lambda^n x^n-1=0
\end{equation}
and the equation (\ref{e7a}) also has exactly one positive root, i.e., $(x,x)$ is a multiple root for the
system of equations (\ref{e5}). This means that, AGM coincide with TIGM.

In particular, if $\lambda x=1$ for $m+r=k-2$ $(n=2)$ or if
$3\lambda^2x^2+2\lambda^3x^3=1$ for $m+r=k-3$ $(n=3)$, then in both cases there is no AGM (not TI).\\

Let
\begin{equation}\label{e8}
\left\{
\begin{array}{ll}
    x=f(y), \\
    y=f(x), \\
    \end{array}
\right.
\end{equation}
where $f(x)=\frac{1}{(1+\lambda x)^{k}}$.

The next lemma is obvious.

\textbf{Lemma 1.} If $(x_0,y_0)$ is a solution to the system of equations (\ref{e8}), then $(y_0,x_0)$ is also a solution to the system of equations (\ref{e8}).

\textbf{Remark 3.} If the solution $(x,y)$ of the system of equations (\ref{e8}) corresponds to alternative Gibbs measure denoted by $\mu$, then the solution $(y,x)$ corresponds to alternative Gibbs measure denoted by $\mu^{'}$.

\subsection{Alternative Gibbs measures in the case $m+r\leq k-2$}\

In this section, we consider the cases $k=2$, $k=3$ and $k=4$. In the case $k=2$ we have only the case $m=0$ and $r=0$. In the case $k=3$ we have $m=0$ and $r=0$; $m=0$ and $r=1$; $m=1$ and $r=0$. In the case $k=4$ we have $m=0$ and $r=0$; $m=0$ and $r=1$; $m=1$ and $r=0$; $m=1$ and $r=1$;
$m=0$ and $r=2$; $m=2$ and $r=0$. In all cases, by Remark 2, we will not consider the case $m=0$ and $r=0$.
corresponds to the translation-invariant Gibbs measure and solutions $(x_1,y_1), (x_2,y_2)$ in Statement 1 ($(x^*,y^*), (y^*,x^*)$ in Statement 2) correspond to two-periodic Gibbs measures (see \protect{\cite{RKhM}}).

\textbf{The case $k=3$, $m=1$ and $r=0$.} For $m=1$ and $r=0$ (resp. $m=0$ and $r=1$) the system of equations (\ref{e5}) can be rewritten
\begin{equation}\label{e11}
\left\{
 \begin{array}{ll}
    x= \frac{1}{1+\lambda x}\cdot\frac{1}{(1+\lambda y)^2}, \\
    y= \frac{1}{(1+\lambda x)^3}.  \\
    \end{array}
\right.
\end{equation}
From the system of equations (\ref{e11}) due to (\ref{e6}) we obtain $(x-y)\Big(\lambda^2xy-1\Big)=0$.
Hence, $x=y$ or $\lambda^2xy=1$. The case $x=y$ has already been considered.

Let $\lambda^2xy=1$. Then $\lambda x=\frac{1}{\lambda y}$ for $x\neq y$. From here and from (\ref{e11}) after some algebras we can get:

\begin{equation}\label{e12}
\left\{
 \begin{array}{ll}
    (1+\lambda x)^3-\lambda^2x=0, \\
    (1+\lambda y)^3-\lambda^3y^2=0.
    \end{array}
\right.
\end{equation}
From $\lambda^2xy=1$ we find $y$ and substitute into the second equation of the system (\ref{e12}). Then
$$
\left\{
 \begin{array}{ll}
    (1+\lambda x)^3-\lambda^2x=0, \\
    \frac{(1+\lambda x)^3-\lambda^2x}{\lambda^3x^3}=0.  \\
   \end{array}
\right.
$$
We introduce the notation $f(x)=(1+\lambda x)^3-\lambda^2x$. Then the roots of the equation $f(x)=0$ are also roots of the system (\ref {e11}).
Using the Cardano formulas, we find the positive solution of the last equation
$$\lambda^3 x^3+3\lambda^2x^2+\Big(3-\lambda^2\Big)x+1=0.$$
Let $x=q-\frac{1}{\lambda}$, then
$$f(q)=\lambda^3q^3-\lambda^2q+\lambda, \ \ D=\frac{1}{\lambda^4}\left(\frac{1}{4}-\frac{1}{27\lambda}\right).$$

If $D>0$, i.e., $\lambda<\frac{27}4$ then by Cardano's formula the equation $f(q)=0$ has one negative root.

If $D=0$,  i.e., $\lambda=\frac{27}4$ then the equation $f(q)=0$ has one multiple positive root of the form $q'=\frac29$, i.e. $x'=\frac{2}{27}$, $y'=\frac{8}{27}$.

By Cardano's formula, the equation $f(q)=0$ has three real roots if $D<0$.
Hence, $f(x)=0$ has three real roots if $\lambda>\frac{27}{4}$.
Let these solutions be $x_1, x_2, x_3$. By the Vieta's formulas we have
$$x_1+x_2+x_3=-\frac{3}{\lambda}, \ \ x_1x_2+x_1x_3+x_2x_3=\frac{3-\lambda^2}{\lambda^3}<0, \ \ x_1x_2x_3=-\frac{1}{\lambda^3}.$$
From equality $x_1x_2+x_1x_3+x_2x_3=\frac{3-\lambda^2}{\lambda^3}$ we obtain that at least one and at most two roots of the equation are positive.

From the equality  $x_1x_2x_3=-1$ it follows that exactly two roots are positive.
Hence, $f(x)=0$ has two positive roots if $\lambda>\frac{27}{4}$.
These roots have the following form
$$x_1=\frac{\sqrt[3]{t^2}-6\sqrt[3]{t}+12\lambda}{6\lambda \sqrt[3]{t}}, \ \
x_2=\frac{6\sqrt[3]{p}}{\lambda(\sqrt[3]{p^2}+(2\lambda-6)\sqrt[3]{p}+4\lambda^2-24\lambda)}.$$
Here
$$t=-108\lambda+12\lambda\sqrt{-12\lambda+81}, \ \ p=108\lambda+8\lambda^3-72\lambda^2+12\lambda\sqrt{-12\lambda+81}.$$

From the equality $\lambda^2xy=1$ we find $y_1$ and $y_2$ corresponding to $x_1$ and $x_2$:
$$y_1=\frac{6\sqrt[3]{t}}{\lambda(\sqrt[3]{t^2}-6\sqrt[3]{t}+12\lambda)}, \ \
y_2=\frac{\sqrt{p^2}+(2\lambda-6)\sqrt[3]{p}+4\lambda^2-24\lambda}{6\sqrt[3]{p}\lambda}.$$

Thus, the following statement is true.

\textbf{Statement 2.} Let $k=3$ and $\lambda_{cr}=\frac{27}{4}$. Then the system of equations (\ref{e11}):

1. for $0<\lambda<\lambda_{cr}$ has a unique solution $(x,x)$;

2. for $\lambda=\lambda_{cr}$ has two solutions $(x,x), (\frac{2}{27},\frac{8}{27})$;

3. for $\lambda>\lambda_{cr}$ has three solutions $(x,x), (x_1,y_1), (x_2,y_2)$.

\textbf{Remark 4.} The measure corresponding to the solution $(x, x)$ is translation invariant and measures corresponding to solutions $(\frac{2}{27},\frac{8}{27})$, $(x_1,y_1), (x_2,y_2)$ are AGMs (not TI).\\

\textbf{Theorem 3.} Let $k=3$ and $r+m\leq 1$, i.e., $m=1$ and $r=0$ or $m=0$ and $r=1$. Then for the HC-model there exists $\lambda_{cr}=\frac{27}{4}$ such that for $0<\lambda<\lambda_{cr}$ there is a unique AGM which coincides with the only TIGM $\mu_0$, for $\lambda=\lambda_{cr}$  there are exactly two AGMs $\mu_0$ and $\mu'$, where $\mu'$ is AGM (not TI) and for $\lambda>\lambda_{cr} $ there are exactly three AGMs $\mu_0$, $\mu_1$ and $\mu_2$, where $\mu_1$ and $\mu_2$ are AGMs (not TI).\\

\textbf{The case $k=4$, $m=1$ and $r=0$ ($m=0$ and $r=1$).} In this case from (\ref{e5}) we get
\begin{equation}\label{e13}
\left\{
 \begin{array}{ll}
    x= \frac{1}{1+\lambda x}\cdot\frac{1}{(1+\lambda y)^3}, \\
    y= \frac{1}{(1+\lambda x)^4}.  \\
    \end{array}
\right.
\end{equation}
From the system of equations (\ref{e13}) due to (\ref{e6}) we can get
$$(y-x)\Big(\lambda^2xy(3+\lambda(x+y))-1\Big)=0.$$
Hence $x=y$ or $\lambda^2xy(3+\lambda(x+y))=1$. It is clear that in the case $x=y$ we obtain a solution corresponding to the TIGM.

Suppose $x\neq y$ and $\lambda^2xy(3+\lambda(x+y))=1$.
Then, substituting the expression for $y$ from the second equation of the system (\ref{e13}) into the last equality, we obtain the equation
$$f(x,\lambda)=\lambda^8x^8+8\lambda^7x^7-\lambda^7x^6+28\lambda^6x^6-7\lambda^6x^5+
56\lambda^5x^5-18\lambda^5x^4+70\lambda^4x^4-22\lambda^4x^3+$$
$$+56\lambda^3x^3-13\lambda^3x^2-\lambda^3x+28\lambda^2x^2-3\lambda^2x+8\lambda x+1=0.$$
Denoting $\lambda x=u$, $u>0$ we then have the equation
$$f(u)=u^8+8u^7+(28-\lambda)u^6+(56-7\lambda)u^5+(70-18\lambda)u^4+(56-22\lambda)u^3+
(28-13\lambda)u^2-(\lambda^2+3\lambda-8)u+1=0,$$
which has a solution $u=u(\lambda)$. But we regard this as an equation for $\lambda$ and obtain  solutions $\lambda=\lambda(u)$:
$$\lambda_1(u)=\frac{(u+1)^4}{2u}\cdot\left(\sqrt{u^4+6u^3+9u^2+4u}-u^2-3u\right),$$
$$\lambda_2(u)=-\frac{(u+1)^4}{2u}\cdot\left(\sqrt{u^4+6u^3+9u^2+4u}+u^2+3u\right).$$
Therefore, because $\lambda_2<0$ for $u>0$, we have
$$\lambda-\lambda_1=0 \ \Rightarrow \ \lambda=\frac{(u+1)^4}{2u}\cdot\left(\sqrt{u^4+6u^3+9u^2+4u}-u^2-3u\right)=\psi(u).$$
Analysis of the function $\psi(u)$ shows that $\psi(u)>0.$
In addition, $\psi(u)\rightarrow +\infty$ as $u\rightarrow
0$ and as $u\rightarrow+\infty$, and each value of $\lambda$ therefore corresponds to at least two values of $u$ for $\lambda>\psi(u^{*})$ but to one value at $\lambda=\psi(u^{*})$, and the equation $\lambda=\psi(u)$ has no solutions for $\lambda<\psi(u^{*}),$ where $u^{*}$ a solution of the equation $\psi'(u)=0$ (see Fig.2). We calculate the derivative
$$\psi'(u)=\frac{(u+1)^3\big[-(5u^2+13u)\sqrt{u^2+4u}+5u^3+23u^2+16u-2\big]}{2\sqrt{u^2+4u}}.$$
It is clear that if $5u^3+23u^2+16u-2<0$ then $\psi^{'}(u)<0$ and the equation $\psi^{'}(u)=0$ has no solutions. So it must be $5u^3+23u^2+16u-2>0$.
$$5u^3+23u^2+16u-2=5(u+1)\left(u+\frac{9-\sqrt{91}}{5}\right)\left(u+\frac{9+\sqrt{91}}{5}\right) \ \ \Rightarrow \ \ u>\frac{\sqrt{91}-9}{5}.$$
We solve the equation $\psi^{'}(u)=0$ for $u>0$:
$$-(5u^2+13u)\sqrt{u^2+4u}+5u^3+23u^2+16u-2=0 \ \ \Rightarrow \ 10u^3+41u^2-16u+1=0.$$
We solve the last equation by the Cardano method:
$$u_1=\frac{\sqrt{2161}}{15}\cdot\cos\left(\frac{\arccos\left(\frac{-99791}{\sqrt{10091699281}}\right)}{3}\right)-\frac{41}{30}\approx0.284824838,$$
$$u_2=\frac{\sqrt{2161}}{15}\cdot\cos\left(\frac{\arccos\left(\frac{-99791}{\sqrt{10091699281}}\right)+2 \pi}{3}\right)-\frac{41}{30}\approx-4.463483795,$$
$$u_3=\frac{\sqrt{2161}}{15}\cdot\cos\left(\frac{\arccos\left(\frac{-99791}{\sqrt{10091699281}}\right)+4 \pi}{3}\right)-\frac{41}{30}\approx0.078658955.$$
Hence, since $u>\frac{\sqrt{91}-9}{5}$ we get the solution $u^{*}=u_1$. We set
$$\lambda_{cr}=\psi(u^{*})\approx2.31.$$

\begin{figure}
\begin{center}
 \includegraphics[width=6cm]{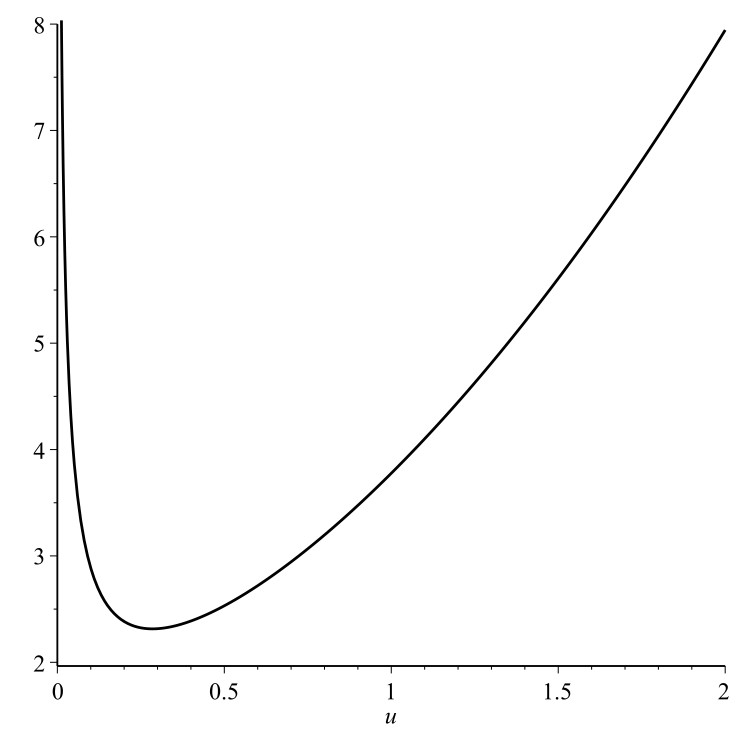}
\end{center}
\begin{center} {\footnotesize \noindent
Figure 2. Graph of the function $\lambda_1(u)$} \
\end{center}
\end{figure}

We note that if $\psi''(u)>0$, then each value of $\lambda$ corresponds to only two values of $u$ for $\lambda>\lambda_{cr}$. We therefore prove that $\psi''(u)>0$. Indeed,
$$\psi''(u)=\frac{2h(u)}{u(u+1)\sqrt{(u^2+4u)^3}}.$$
Here
$$h(u)=5u^8+56u^7+234u^6+463u^5+460u^4+210u^3+26u^2-u+3-(5u^2+11u)\sqrt{(u^4+6u^3+9u^2+4u)^3}.$$
From the inequality $h(u)>0$ for $u>0$ we obtain
$$\big(5u^8+56u^7+234u^6+463u^5+460u^4+210u^3+26u^2-u+3\big)^2-
\big(5u^2+11u\big)^2\big(u^4+6u^3+9u^2+4u\big)^3=$$
$$=10u^{13}+182u^{12}+1372u^{11}+5505u^{10}+12786u^9+17913u^8+15564u^7+$$
$$+9186u^6+5034u^5+3016u^4+1208u^3+156u^2+(u-3)^2>0.$$
Thus, each value of $\lambda$ corresponds to only two values of $u$ for $\lambda>\lambda_{cr}$.

This can also be seen by computer analysis, i.e., computer analysis shows that the equation $f(x,\lambda)=0$ for $\lambda<\lambda_{cr}$ has no positive solution, at
$\lambda=\lambda_{cr}$ has one positive solution and for $\lambda>\lambda_{cr}$ there are exactly two positive solutions (see Fig. 3).
\begin{figure}
\begin{center}
 \includegraphics[width=6cm]{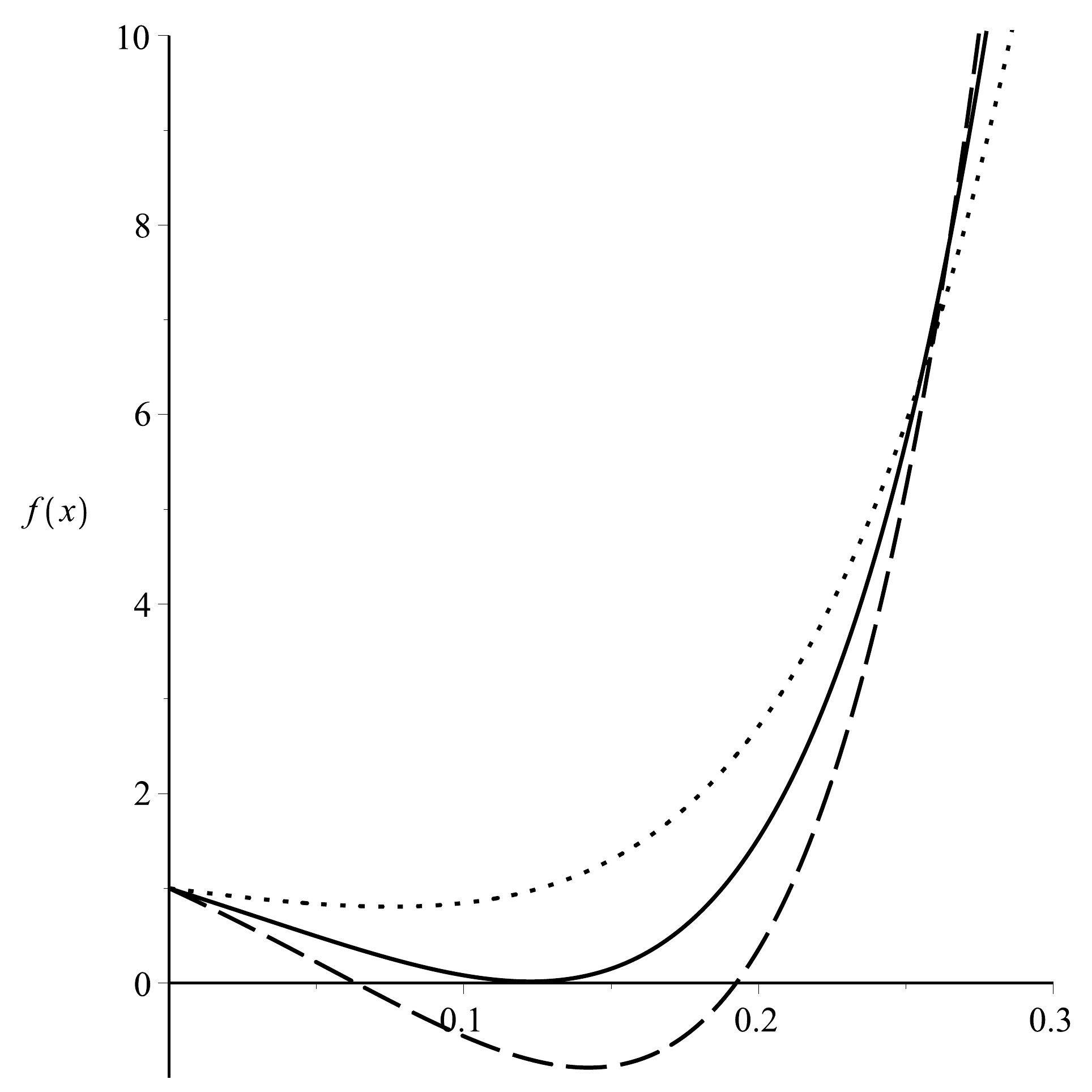}
\end{center}
\begin{center} {\footnotesize \noindent
Figure 3. Graph of the function $f(x,2)$ (dotted line), $f(x,2.3143)$ (continuous line) and $f(x,2.5)$ (dashed line).} \
\end{center}
\end{figure}

Thus, the following statement is true.

\textbf{Statement 3.} Let $k=4$ and $\lambda_{cr}\approx 2.31$. Then the system of equations (\ref{e13}):

1. for $0<\lambda<\lambda_{cr}$ has a unique solution $(x,x)$;

2. for $\lambda=\lambda_{cr}$ has two solutions $(x,x), (x',y')$;

3. for $\lambda>\lambda_{cr}$ has three solutions $(x,x), (x_1,y_1), (x_2,y_2)$.

\textbf{Remark 5.} The measures corresponding to the solution in the Statement 3 for $x\neq y$ are AGMs (not periodic) and they different from previous AGMs.\\

\textbf{The case $k=4$, $m=1$ and $r=1$.} In this case from the system of equations (\ref{e5}) we obtain
\begin{equation}\label{e14}
\left\{
 \begin{array}{ll}
    x= \frac{1}{1+\lambda x}\cdot\frac{1}{(1+\lambda y)^3}, \\
    y= \frac{1}{1+\lambda y}\cdot\frac{1}{(1+\lambda x)^3}.  \\
    \end{array}
\right.
\end{equation}
From (\ref{e14}) due to (\ref {e6}) we can get
$$(x-y)\Big(\lambda^2xy-1\Big)=0.$$
Hence $x=y$ or $\lambda^2xy=1$. The case $x=y$ corresponds to the only TIGM.

Let $x\neq y$ and $\lambda^2xy=1$, i.e., $\lambda x=\frac{1}{\lambda y}$. After some algebras the system of equations (\ref{e14}) has the form
\begin{equation}\label{e15}
\left\{
 \begin{array}{ll}
    (1+\lambda x)^4-\lambda^3x^2=0, \\
    (1+\lambda y)^4-\lambda^3y^2=0.
 \end{array}
\right.
\end{equation}

Obviously, that the roots of the equation $f(x)=(1+\lambda x)^4-\lambda^3x^2=0 $ are also roots of (\ref{e14}).
The solutions of the equations $f(x)=0$ and $f(y)=0$ have the form
$$x_{1,2}=\frac{\sqrt{\lambda}-2\pm\sqrt{\lambda-4\sqrt{\lambda}}}{2\lambda}, \ \ y_{1,2}=\frac{\sqrt{\lambda}-2\pm\sqrt{\lambda-4\sqrt{\lambda}}}{2\lambda}.$$
It is easy to see that $x_{1,2}>0$ ($y_{1,2}>0$) for $\lambda\geq16$, and they take complex values for $\lambda<16$. Moreover, $x_1=x_2$ ($y_1=y_2$) for $\lambda=16$ and it coincides with the only translation-invariant solution of (\ref{e14}).

By virtue of the equation $\lambda^2xy=1$ and Lemma 1, we obtain that in the case $x\neq y$ the system of equations (\ref{e14}) has solutions of the form $(x,y)$ and $(y,x)$ for $\lambda>\lambda_{cr}=16$, where
$$x=x_{1}=\frac{\sqrt{\lambda}-2+\sqrt{\lambda-4\sqrt{\lambda}}}{2\lambda}, \ \ y=y_{2}=\frac{\sqrt{\lambda}-2-\sqrt{\lambda-4\sqrt{\lambda}}}{2\lambda},$$
$$y=x_{2}=\frac{\sqrt{\lambda}-2-\sqrt{\lambda-4\sqrt{\lambda}}}{2\lambda}, \ \ x=y_{1}=\frac{\sqrt{\lambda}-2+\sqrt{\lambda-4\sqrt{\lambda}}}{2\lambda}.$$
Thus, the following statement holds

\textbf{Statement 4.} Let $k=4$ and $\lambda_{cr}=16$. Then the system of equations (\ref{e14}):

1. for $0<\lambda\leq\lambda_{cr}$ has a unique solution $(x,x)$;

2. for $\lambda>\lambda_{cr}$ has three solutions $(x,x), (x,y), (y,x)$.

\textbf{Remark 6.} The measure corresponding to the solution $(x,y), (y,x)$ in the Statement 4 are AGMs (not periodic) and they different from previous AGMs.\\

\textbf{The case $k=4$, $m=2$ and $r=0$ ($m=0$ and $r=2$).} In this case from (\ref{e5}) we obtain
\begin{equation}\label{e16}
\left\{
 \begin{array}{ll}
    x= \frac{1}{(1+\lambda x)^2}\cdot\frac{1}{(1+\lambda y)^2}, \\
    y= \frac{1}{(1+\lambda x)^4}.  \\
    \end{array}
\right.
\end{equation}

Using (\ref{e6}) from (\ref{e16}) we can get
$$(x-y)\Big(\lambda^2xy-1\Big)=0.$$
Hence, $x=y$ or $\lambda^2xy=1$. The case $x=y$ corresponds to the only TIGM.

We consider the case $x\neq y$ and $\lambda^2xy=1$ $\Big(\lambda x=\frac{1}{\lambda y}\Big)$. After some algebras (\ref{e16}) has the form
\begin{equation} \label{e17}
\left\{
 \begin{array}{ll}
    (1+\lambda x)^4-\lambda^2x=0, \\
    (1+\lambda y)^4-\lambda^4y^3=0.
 \end{array}
\right.
\end{equation}
From the equation $\lambda^2xy=1$ we find $y$ and substitute it for the second equation (\ref{e17}). Then
$$
\left\{
 \begin{array}{ll}
    (1+\lambda x)^4-\lambda^2x=0, \\
    \frac{(1+\lambda x)^4-\lambda^2x}{\lambda^4x^4}=0.  \\
   \end{array}
\right.
$$
Let's rewrite the equation $f(x)=(1+\lambda x)^4-\lambda^2x=0$ as
$$\lambda^4x^4+4\lambda^3x^3+6\lambda^2x^2 +\lambda(4-\lambda)x+1=0.$$
We solve the last equation by the Ferrari method from linear algebra. We introduce the notation $x=t-\frac{1}{\lambda}$. Then
$$f\left(t-\frac{1}{\lambda}\right)=\lambda^4t^4-\lambda^2t+\lambda=
(\lambda^2t^2+p)^2-2\lambda^2p\Big(t+\frac{1}{4p}\Big)^2=$$ $$=\left(\lambda^2t^2+p-\lambda \sqrt {2p}\Big(t+\frac{1}{4p}\Big)\right) \left(\lambda^2t^2+p+\lambda
\sqrt {2p}\Big(t+\frac{1}{4p}\Big)\right)=0,$$
where $$p=\frac{\sqrt[3]{108\lambda^2+12\sqrt{81\lambda^4-768\lambda^3}}}{12}+
\frac{4\lambda}{\sqrt[3]{108\lambda^2+12\sqrt{81\lambda^4-768\lambda^3}}}.$$
Solutions have the following form
$$t_{1,2}=\frac{\sqrt{2p^3}\pm\sqrt{\sqrt{2p^3}\lambda-2p^3}}{2\lambda p}, \ \  t_{3,4}=\frac{-\sqrt{2p^3}\pm\sqrt{-\sqrt{2p^3}\lambda-2p^3}}{2\lambda p}.$$
By virtue $x=t-\frac{1}{\lambda}$, for solutions we obtain
$$x_{1,2}=\frac{\sqrt{2p^3}\pm\sqrt{\sqrt{2p^4}\lambda-2p^3}-2p}{2\lambda p}, \ \  x_{3,4}=\frac{-\sqrt{2p^3}\pm\sqrt{-\sqrt{2p^3}\lambda-2p^3}-2p}{2\lambda p}.$$
Computer analysis shows that $x_{1,2}>0$ for $\lambda>\lambda_{cr}\approx9.48$, and values $x_{3,4}$ are negative or take on complex values for $\lambda>0$ (see Fig. 4).
\begin{figure}
\begin{center}
 \includegraphics[width=6cm]{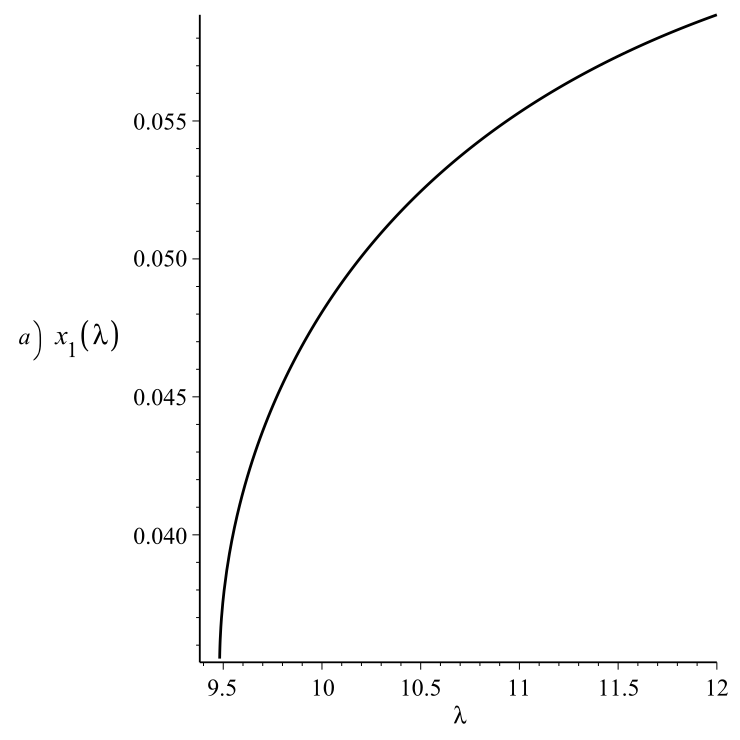} \ \ \ \includegraphics[width=6cm]{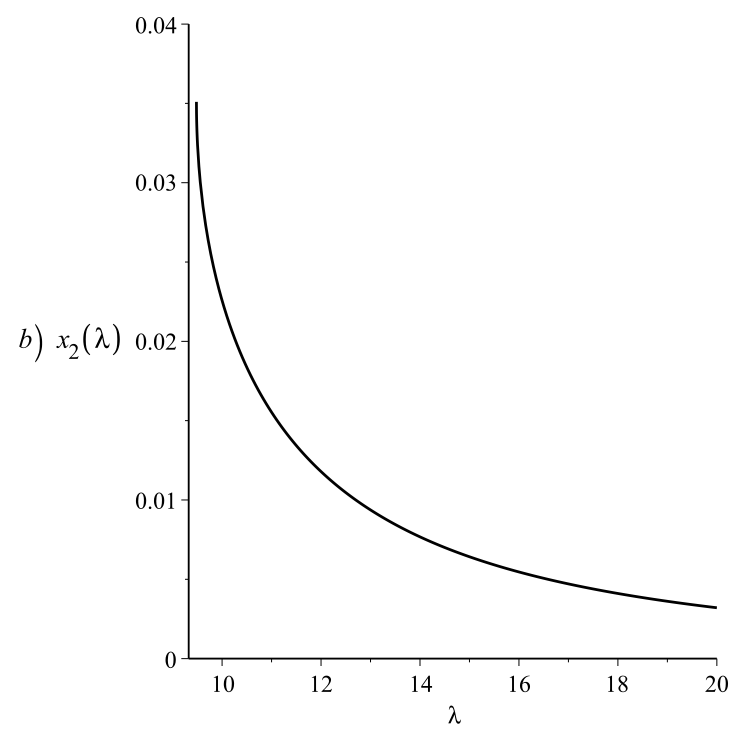}
\end{center}
\begin{center}{\footnotesize \noindent
Figure 4. a) Graph of the function $x_1(\lambda)$ at $\lambda\in[9.4;12]$, b) Graph of the function $x_2(\lambda)$ at $\lambda\in[9.4;20]$.}\
\end{center}
\end{figure}
Values $y_1$ and $y_2$ corresponding to values $x_1$ and $x_2$ have the form:
$$y_{1,2}=\frac{2p}{\lambda \Big(\sqrt{2p^3}\pm\sqrt{\sqrt{2p^3}\lambda-2p^3}-2p\Big)}.$$
For $\lambda=\lambda_{cr}\approx9.4815$ the system of equations (\ref{e17}) has solutions of the form
$$x'=\frac{\sqrt{2p^3}-2p}{2\lambda p}, \ \ y'=\frac{2p}{\lambda \sqrt{2p^3}-2p}.$$

Thus, the following statement is true.

\textbf{Statement 5.} Let $k=4$ and $\lambda_{cr}\approx9.48$. Then the system of equations (\ref{e17}):

1. for $0<\lambda<\lambda_{cr}$ has a unique solution $(x,x)$;

2. for $\lambda=\lambda_{cr}$ has two solutions $(x,x), (x',y')$;

3. for $\lambda>\lambda_{cr}$ has three solutions $(x,x), (x_1,y_1), (x_2,y_2)$.

By using all propositions, we get the following theorem.

\textbf{Theorem 4.} Let $k=4$ and $r+m\leq 2$. For the HC-model the following statements are true:

1. If $m=1$ and $r=0$ or $m=0$ and $r=1$, then  there exists $\lambda_{cr}\approx2.31$ such that for $0<\lambda<\lambda_{cr}$ there is a unique AGM which coincides with the only TIGM $\mu_0$,  for $\lambda=\lambda_{cr}$  there are exactly two AGMs $\mu_0$ and $\mu'$, where $\mu'$ is AGM (not TI) and for $\lambda>\lambda_{cr} $ there are exactly three AGMs $\mu_0$, $\mu_1$ and $\mu_2$, where $\mu_1$ and $\mu_2$ are AGMs (not TI).

2. If $m=1$ and $r=1$ then there exists $\lambda_{cr}=16$ such that for $0<\lambda\leq\lambda_{cr}$ there is a unique AGM which coincides with the only TIGM $\mu_0$, for $\lambda>\lambda_{cr} $ there are exactly three AGMs $\mu_0$, $\mu_1$ and $\mu_2$, where $\mu_1$ and $\mu_2$ are AGMs (not TI).

3. If $m=2$ and $r=0$ or $m=0$ and $r=2$, then there exists $\lambda_{cr}\approx9.48$ such that for $0<\lambda<\lambda_{cr}$ there is a unique AGM which coincides with the only TIGM $\mu_0$, for $\lambda=\lambda_{cr}$  there are exactly two AGMs $\mu_0$ and $\mu'$, where $\mu'$ is AGM (not TI) and for $\lambda>\lambda_{cr} $ there are exactly three AGMs $\mu_0$, $\mu_1$ and $\mu_2$, where $\mu_1$ and $\mu_2$ are AGMs (not TI).\\

\section{The case $m+r\leq k-2$ ($m=r$)}\

The following lemma is known.

\textbf{Lemma 2.} \protect{\cite{K}} \textit{Let $f:[0,1]\rightarrow [0,1]$
be a continuous function with a fixed point $\xi \in (0,1)$. We
assume that $f$ is differentiable at $\xi$ and $f^{'}(\xi)<-1.$
Then there exist points $x_0$ and $x_1$, $ 0\leq x_0<\xi<x_1
\leq1,$ such that $f(x_0)=x_1$ and $f(x_1)=x_0.$}\

For $m=r$ by (\ref{e5}) we obtain
\begin{equation}\label{eq8}
\begin{cases}
    x=\frac{1}{(1+\lambda x)^m}\cdot \frac{1}{(1+\lambda y)^{k-m}}; \\
    y=\frac{1}{(1+\lambda y)^{m}}\cdot \frac{1}{(1+\lambda x)^{k-m}}.
    \end{cases}
\end{equation}
Here $x,y\in(0;1)$.
After some transformations from (\ref{eq8}) we obtain the following system of equations:
\begin{equation}\label{eq9}
\begin{cases}
    y=f(x); \\
    x=f(y),
\end{cases}
\end{equation}
where
$$f(x)=\frac1\lambda\cdot\left(\frac{1}{x(1+\lambda x)^{m}}\right)^{\frac1{k-m}}-\frac1\lambda.$$
From (\ref{eq9}) we get the equation $f(f(x))=x$.

First, we consider the equation $f(x)=x$. The function $f(x)$ is differentiable and decreasing for $0<x<1$:
$$f'(x)=-\frac{1+\lambda(m+1)x}{\lambda (k-m)x^\frac{k-m+1}{k-m}(1+\lambda x)^\frac{k}{k-m}}<0.$$

We rewrite the equation $f(x)=x$:
$$x=\frac1\lambda\cdot\left(\frac{1}{x(1+\lambda x)^{m}}\right)^{\frac1{k-m}}-\frac1\lambda\ \ \Rightarrow \ \  (1+\lambda x)^k=\frac1x.$$

It is known from \protect{\cite{7}} that the last equation has a unique solution $\tilde{x}$, i.e., the equation $f(x)=x$ has a unique solution $\tilde{x}$.

We solve the inequality $f'(\tilde{x})<-1$:
$$-\frac{1+\lambda(m+1)\tilde{x}}{\lambda (k-m)\tilde{x}^\frac{k-m+1}{k-m}(1+\lambda \tilde{x})^\frac{k}{k-m}}<-1 \ \ \Rightarrow \ \
\frac{1+(m+1)\lambda \tilde{x}}{\lambda(k-m)\tilde{x}}>1 \ \ \Rightarrow \ \tilde{x}<\frac{1}{\lambda(k-2m-1)}.$$
Then from $f(\tilde{x})=\tilde{x}$ we get
$$\left(1+\frac1{k-2m-1}\right)^k<\lambda(k-2m-1) \ \ \Rightarrow \ \ \lambda>\lambda_{cr}=\left(\frac{k-2m}{k-2m-1}\right)^k\cdot\frac1{k-2m-1}.$$
Hence, by Lemma 1 and Lemma 2 the system of equations (\ref{eq8}) for $\lambda>\lambda_{cr}$ has at least three positive solutions $(x, y), \ (\tilde{x}, \tilde{x}), \ (y, x)$, where $x\neq y$.

Thus, the following theorem is true.

\textbf{Theorem 5.} Let $k\geq2$, $m+r\leq k-2$ ($m=r$) and $\lambda_{cr}=\left(\frac{k-2m}{k-2m-1}\right)^k\cdot\frac1{k-2m-1}$. Then for the HC-model for $\lambda>\lambda_{cr}$ there are at least three Gibbs measures one of which is TI and the other are AGMs (not TI).

\subsection{The case $m+r=k-2$, $k\geq2$}\

In the case $m+r=k-2$ $(n=2)$, the system of equations (\ref{e5}) has the form:
\begin{equation}\label{e18}
\left\{
 \begin{array}{ll}
    x=\frac{1}{(1+\lambda x)^m}\cdot \frac{1}{(1+\lambda y)^{k-m}}; \\
    y=\frac{1}{(1+\lambda y)^{k-m-2}}\cdot \frac{1}{(1+\lambda x)^{m+2}}.  \\
    \end{array}
\right.
\end{equation}
From the system of equations (\ref{e18}) due to (\ref{e6}) we can get
$$(x-y)\Big(\lambda^2xy-1\Big)=0.$$
Hence, $x=y$ or $\lambda^2xy=1$. The case $x=y$ has already been considered.

Let $\lambda^2xy=1$. Then $\lambda x=\frac{1}{\lambda y}$ for $x\neq y$. By virtue (\ref{e18}), after some algebras, we can obtain the system of equations
$$
\left\{
 \begin{array}{ll}
    x=\frac{(\lambda x)^{k-m}}{(1+\lambda x)^{k}}, \\[2mm]
    y=\frac{(\lambda y)^{m+2}}{(1+\lambda y)^{k}},
    \end{array}
\right.
$$
which is equivalent to the system of equations:
\begin{equation}\label{e19}\left\{
 \begin{array}{ll}
    (1+\lambda x)^{k}-\lambda^{k-m}x^{k-m-1}=0, \\
    (1+\lambda y)^{k}-\lambda^{m+2}y^{m+1}=0.
    \end{array}
\right.
\end{equation}
From the equation $\lambda^2xy=1$ we find $y$ and substitute it for the second equation of the system equations (\ref{e19}). Then
$$
\left\{
 \begin{array}{ll}
    (1+\lambda x)^{k}-\lambda^{k-m}x^{k-m-1}=0, \\ [2mm]
    \frac{(1+\lambda x)^{k}-\lambda^{k-m}x^{k-m-1}}{\lambda^kx^k}=0.
   \end{array}
\right.
$$
We consider the function
$$f(x)=(1+\lambda x)^{k}-\lambda^{k-m}x^{k-m-1}.$$
Obviously, the roots of the equation $f(x)=0$ are also roots of (\ref{e19}). Let's rewrite $f(x)$ as a polynomial:
$$f(x)=\lambda^kx^k+C_k^1\lambda^{k-1}x^{k-1}+\cdots+C_k^{m+1}\lambda^{k-m-1}x^{k-m-1}+\cdots+C_k^{k-1}\lambda x+1-\lambda^{k-m}x^{k-m-1}$$
or
$$f(x)=\lambda^kx^k+C_k^1\lambda^{k-1}x^{k-1}+\cdots+(C_k^{m+1}-\lambda)\lambda^{k-m-1}x^{k-m-1}+\cdots+C_k^{k-1}\lambda x+1.$$
If $\lambda<C_k^{m+1}$, then $f(x)=0$ has no positive solutions, if $\lambda>C_k^{m+1}$, then number of sign changes of the first equation of the last equality is
two. Due to the Descartes' theorem, the equation $f(x)=0$ has at most two positive solutions.

On the other hand, it is easy to see that $0<x<1$, $f(0)=1$ and $f(1)=(1+\lambda)^k-\lambda^{k-m}>0.$  Moreover, $f\Big(\frac{1}{\lambda}\Big)=2^k-\lambda<0$, if
$\lambda>2^k.$

It follows from the above that there exists $\lambda_{cr}: C_k^{m+1}<\lambda_{cr}\leq2^k $ such that for $\lambda>\lambda_{cr}$ the equation $f(x)=0$ has two positive solutions, for $\lambda=\lambda_{cr}$ has twice multiplicity positive solution and $\lambda<\lambda_{cr}$  has no positive
solution.

When $m=r$, the system of equations (\ref{e18}) can be written as
\begin{equation}\label{e20}
\left\{
 \begin{array}{ll}
    x=\frac{1}{(1+\lambda x)^m}\cdot \frac{1}{(1+\lambda y)^{k-m}}; \\
    y=\frac{1}{(1+\lambda y)^{m}}\cdot \frac{1}{(1+\lambda x)^{k-m}}.  \\
    \end{array}
\right.
\end{equation}
It follows from the Lemma 1 that if the number of solutions of the equation $x=f(x)$ is odd or even, then the number of solutions of $x=f(f(x))$ is also respectively odd or even.

As a result, when $m=r$ and $\lambda=\lambda_{cr}$, the number of solutions of the system of equations (\ref{e20}) cannot be even. Because the TI solution was
unique. It follows that the solution of the system of equations (\ref{e20}) corresponding to  $\lambda=\lambda_{cr}$  coincides  with the translation-invariant solution.

Thus, we have proved the following theorem.

\textbf{Theorem 6.} Let $k\geq2$ and $r+m=k-2$. Then for the HC-model there exists $\lambda_{cr}$ such that next statements are true:

1. for $0<\lambda<\lambda_{cr}$ there is a unique AGM and it coincides with the only TIGM $\mu_0$;

2. if $m=r$ and $\lambda=\lambda_{cr}$ then there is a unique AGM and it coincides with the only TIGM $\mu_0$;

3. if $m\neq r$ and $\lambda=\lambda_{cr}$ there is at least one AGM;

4. for $\lambda>\lambda_{cr} $ there are exactly three Gibbs measures $\mu_0$, $\mu_1$ and $\mu_2$, where $\mu_1$ and $\mu_2$ are AGMs (not TI).\\

\section{Relation of the Alternative Gibbs measures to known ones}\

\textit{Translation invariant measures.} (see \protect{\cite{7}}) Such measures correspond to $z_x\equiv z$, i.e. constant functions. These measures are particular cases of
our measures mentioned which can be obtained for $m=k$, i.e. $k-m=0$.
In this case the condition (\ref{e3}) reads
\begin{equation}
\label{f1} z=\frac{1}{(1+\lambda z)^{k}}.
\end{equation}
The equation (\ref{f1}) has a unique solution for all $\lambda>0.$

\textit{Bleher-Ganikhodjaev construction}. Consider an infinite path $\pi=\{x^0=x_0<x_1<...\}$ on the half Cayley tree (the notation $x<y$ meaning that paths from
the root to $y$ go through $x$). Associate to this path a collection $z^{\pi}$ of numbers given by the condition
$$\label{f2}
z^{\pi}_x=
\begin{cases} l ~\mbox{if}~ x\prec x_n,~x\in W_n,\\
 h, ~\mbox{if}~ x_n\prec x,~x\in W_n,\\
 h, ~\mbox{if}~ x=x_n.
 \end{cases}
 $$
$n=1,2,...$ where $x\prec x_n$ (resp. $x_n\prec x$) means that $x$ is on the left (resp. right) from the path $\pi$ and $z_{x_n}\in \{h,l\}$ are arbitrary numbers.
For any infinite path $\pi$, the collection of numbers $z^{\pi}$ satisfying relations (\ref{e3}) exists and is unique (see Fig. 5).

\begin{figure}
\scalebox{1}{
\begin{tikzpicture}[level distance=2.5cm,
level 1/.style={sibling distance=2.3cm},
level 2/.style={sibling distance=2.5cm},
level 3/.style={sibling distance=.3cm}]
\node {$h$} [grow'=up]
    child[sibling distance=3.5cm] {node {$h$}
        child[sibling distance=.5cm] {node {$h$}}
        child[sibling distance=.5cm] {node {$h$}}
        child[sibling distance=.5cm] {node {$h$}}
        child[sibling distance=.5cm, very thick] {node {$h$}}
        child[sibling distance=.5cm] {node {$l$}}
    }
    child[sibling distance=3.4cm] {node {$h$}
        child[sibling distance=.5cm] {node {$h$}}
        child[sibling distance=.5cm] {node {$h$}}
        child[sibling distance=.5cm] {node {$h$}}
        child[sibling distance=.5cm, very thick] {node {$h$}}
        child[sibling distance=.5cm] {node {$l$}}
    }
    child[sibling distance=2.5cm] {node {$h$}
        child[sibling distance=.5cm] {node {$h$}}
        child[sibling distance=.5cm] {node {$h$}}
        child[sibling distance=.5cm] {node {$h$}}
        child[sibling distance=.5cm, very thick] {node {$h$}}
        child[sibling distance=.5cm] {node {$l$}}
    }
    child[sibling distance=3.4cm, very thick]{node {$h$}
        child[sibling distance=.5cm, thin] {node {$h$}}
        child[sibling distance=.5cm, thin] {node {$h$}}
        child[sibling distance=.5cm, thin] {node {$h$}}
        child[sibling distance=.5cm, very thick] {node {$h$}}
        child[sibling distance=.5cm, thin] {node {$l$}}
    }
    child[sibling distance=3.5cm] {node {$l$}
        child[sibling distance=.5cm] {node {$h$}}
        child[sibling distance=.5cm] {node {$h$}}
        child[sibling distance=.5cm] {node {$h$}}
        child[sibling distance=.5cm, very thick] {node {$h$}}
        child[sibling distance=.5cm] {node {$l$}}
    }
;
\end{tikzpicture}
}
\begin{center}{\footnotesize \noindent Figure 5. In this figure the values of function $z_x$ on the vertices of the Cayley tree
of order 5 are shown. This is the case when $m=4$ and $r=4$.}
\end{center}
\end{figure}

\textit{Periodic Gibbs measures.} (see \protect{\cite{7}}) Let $G_k$ be a free product of $k+1$ cyclic groups of the second order with generators $a_1,a_2,...,a_{k+1},$
respectively.

It is known that there exists an one-to-one correspondence between the set of vertices $V$ of the Cayley tree $\Im^k$ and the group $G_k$.

\textbf{Definition 6.}  Let $\widehat{G}$ be a normal subgroup of the group $G_k$. The set $z=\{z_x,x\in G_k\}$
is said to be $\widehat{G}$ -periodic if $z_{yx}=z_x$ for $\forall x\in G_k, y\in\widehat{G}_k.$

Let $G^{(2)}_k=\{x\in G_k : \mbox{the length of word} ~x~ \mbox{is even} \}.$
Note that $G^{(2)}_k$ is the set of even vertices (i.e. with even distance to the root). Consider the boundary condition $h$ and $l$:
$$z_x=
\begin{cases} h ~ \mbox{if}~ x\in G^{(2)}_k,\\
 l ~ \mbox{if}~ x\in G_k\setminus G^{(2)}_k.
 \end{cases}$$
and denote by $\mu_{1}$, $\mu_{2}$ the corresponding Gibbs measures.
The $\widehat{G}$- periodic solutions of equation (\ref{e3}) are either translation-invariant ($G_k$- periodic) or $G^{(2)}_k$
-periodic, they are solutions to
$$\begin{cases} h=\frac{1}{(1+\lambda l)^{k}}, \\
l=\frac{1}{(1+\lambda h)^{k}}.
 \end{cases}$$
We note that these measures are particular cases of measures of $\mu_{h,l}$ which
can be obtained for $m=r=0$ (See figure 6, for $k=4$).

\begin{figure}
\scalebox{1}{
\begin{tikzpicture}[level distance=2.5cm,
level 1/.style={sibling distance=2.3cm},
level 2/.style={sibling distance=2.5cm},
level 3/.style={sibling distance=.3cm}]
\node {$h$} [grow'=up]
    child[sibling distance=4.5cm] {node {$l$}
        child[sibling distance=.6cm] foreach \name in {h,h,h,h} { node
        {$\name$} }
    }
    child[sibling distance=4.5cm] {node {$l$}
        child[sibling distance=.6cm] foreach \name in {h,h,h,h} { node
        {$\name$} }
    }
    child[sibling distance=4.5cm] {node {$l$}
        child[sibling distance=.6cm] foreach \name in {h,h,h,h} { node
        {$\name$} }
    }
    child[sibling distance=4.5cm] {node {$l$}
        child[sibling distance=.6cm] foreach \name in {h,h,h,h} { node
        {$\name$} }
    }
;
\end{tikzpicture}
}
\begin{center}{\footnotesize \noindent  Figure 6. In this figure the values of function $z_x$ on the vertices of the Cayley tree
of order 4 are shown.}
\end{center}
\end{figure}

\textit{Weakly periodic Gibbs measures.} Following \protect{\cite{Kh1},\cite{Kh2}, \cite{KhRM}} recall the notion of weakly periodic Gibbs measures. Let
$G_k/\widehat{G}_k=\{H_1, ..., H_r\}$ be a factor group, where $\widehat{G}_k$ is a normal subgroup of index $r>1$.

\textbf{Definition 7.}  A set $z=\{z_x,x\in G_k\}$ is called $\widehat{G}_k$ - weakly periodic, if $z_x=z_{ij}$, for any
$x\in H_i$, $x_\downarrow \in H_j$, where $x_\downarrow$ denotes the ancestor of $x$.

We recall results known for the cases of index two. Note that any such subgroup
has the form
$$H_A=\Big\{x\in G_k:\sum_{i\in A}{w_x(a_i)}~\mbox{is even}\Big\}$$
where ${\emptyset}\neq A\subseteq N_k=\{1,2,...,k+1\}$, and $w_x(a_i)$ is the number of $a_i$ in a word $x\in G_k.$
We consider $A\neq N_k$: when $A=N_k$ weak periodicity coincides with standard periodicity. Let $G_k/H_A=\{H_0,H_1\}$ be the factor group, where $H_0=H_A$,
$H_1=G_k\setminus H_A$. Then, in view of (\ref{e3}), the $H_A$ - weakly periodic b.c. has the form
$$z_x=
\begin{cases} z_1, \ x\in H_A, \ x_\downarrow \in H_A, \\
z_2, \ x\in H_A, \ x_\downarrow \in G_k\setminus H_A,\\
z_3, \ x\in G_k\setminus H_A, \ x_\downarrow \in H_A, \\
z_4, \ x\in G_k\setminus H_A, \ x_\downarrow \in G_k\setminus H_A.
\end{cases}$$
where the $h_i$ satisfy the following equations:
$$z_1=\frac{1}{\Big(1+\lambda z_3\Big)^i}\frac{1}{\Big(1+\lambda z_1\Big)^{k-i}}, \ \ \
z_2=\frac{1}{\Big(1+\lambda z_3\Big)^{i-1}}\frac{1}{\Big(1+\lambda z_1\Big)^{k-i+1}},$$
\begin{equation}\label{f3} z_3=\frac{1}{\Big(1+\lambda z_2\Big)^{i-1}}\frac{1}{\Big(1+\lambda z_4\Big)^{k-i+1}}, \ \
z_4=\frac{1}{\Big(1+\lambda z_2\Big)^{i}}\frac{1}{\Big(1+\lambda z_4\Big)^{k-i}}.\end{equation}
It is obvious that the following sets are invariant with respect to the operator $W:R^4\rightarrow R^4$ defined by RHS of (\ref{f3}):
$$I_1=\Big\{z\in R^4:z_1=z_2=z_3=z_4\Big\}, \ \ \ \ I_2=\Big\{z\in R^4: z_1=z_4; z_2=z_3\Big\}$$
It is obvious to see that

$\bullet$ measures corresponding to solutions on $I_1$ are translation invariant

$\bullet$ measures corresponding to solutions on $I_2$ are weakly periodic, which coincide with the measures given for $m=k-i$, $k-m=i$, $r=i-1$, $k-r=k-i+1$.

\section{Free energy}

 In this section, we consider free energy of \emph{HC}-model Gibbs measure. In fact, Gibbs measures give the probability of the system $X$ being in state $x\in X$ (equivalently, of the random variable $X$ having value $x$) as
$$\mu(X=x)=\frac{1}{Z(\beta)}\exp(-\beta H(x)),$$
where $H(x)$ is a function from the space of states to the real numbers. The parameter $\beta$ is (a free parameter) the inverse temperature. The normalizing constant $Z(\beta)$ is the partition function.

Consider an infinite graph $G$, and let $\Lambda\subset G$ be finite subset. It is convinient to work with reduced free energy $f=-\beta F$, which per unit volume is
$$f(\beta, \Lambda)=\frac{1}{|\Lambda|}\ln Z(\beta, \Lambda),$$
where $Z(\beta, \Lambda)$ is the restiriction of the partition function $Z(\beta)$ on the set $\Lambda$, by fixing the state of the system outside of $\Lambda$.

Note that by Theorem 6, we can construct Aternating Gibbs measures and by using these measures we compute the free energy for such measures.
From \protect{\cite{GRRR, GRHR}} it's known that the free energy of a compatible boundary condition (b.c.) is  defined as the limit:
\begin{equation}\label{fe1}
F(h)=-\lim_{n\to \infty}\frac{1}{ \beta |V_n|}
\ln Z_n
\end{equation}
if it exists.
Here $|\cdot|$ denotes  the cardinality of a set and $Z_n$ is a partition function.
We recall that in our case:
\begin{equation}\label{fe2} Z_n=\sum_{{\widetilde\sigma}_n\in\Omega_{V_n}}
\lambda^{\#{\widetilde\sigma}_n}\prod_{x\in W_n}
z_{{\widetilde\sigma}(x),x}.\end{equation}

We consider ALT Gibbs measures on the half tree and from above, the family of probability measures ce compatible iff $z=\{z_x, x\in G_k\}$ satisfies the equality (3). Also, we considered a special class of $z=\{z_x, x\in G_k\}$ such that:

$\bullet$ if at vertex $x$ we have $z_x=h$, then the function $z_y$, which gives real values to each vertex $y\in S(x)$ by the following rule
$$\begin{cases}
h ~\mbox{on} ~m~ \mbox{vertices of} ~ S(x),\\
l ~\mbox{on} ~k-m ~\mbox{remaining vertices,}
\end{cases}$$

$\bullet$ if at vertex $x$ we have $z_x=l$, then the function $z_y$, which gives real values to each vertex $y\in S(x)$ by the following rule
$$\begin{cases}
l ~\mbox{on} ~r~ \mbox{vertices of} ~ S(x),\\
h ~\mbox{on} ~k-r ~\mbox{remaining vertices.}
\end{cases}$$

 Denote
\begin{equation}\label{ss}
\alpha_n=|\{x\in W_n: z_x=h\}|; \ \ \beta_n=|\{x\in W_n: z_x=l\}|.
\end{equation}
Recall that $W_n$ is the sphere with the center $x^0$ and radius $n$ on the half tree.

Consequently, the following recurrence system holds
\begin{equation}\label{oab}
\left\{\begin{array}{lllll}
\alpha_{n+1}=m\alpha_n+(k-r)\beta_n\\[3mm]
\beta_{n+1}=(k-m)\alpha_n+r\beta_n.
\end{array}\right.
\end{equation}
Denoting $\varphi_n=\alpha_n+\beta_n$, from (\ref{oab})
one gets
\begin{equation}\label{oa1}
\varphi_{n+1}=k\varphi_n \ \Rightarrow \ \varphi_n=k^n, \ n\in \mathbb{N}.
\end{equation}
Since $\alpha_n=k^n-\beta_n$, we get
\begin{equation}\label{eq}k^{n+1}-mk^n=(k-m-r)\beta_n+\beta_{n+1}.\end{equation}
Put
$$\beta_n=\frac{(k-m)k^n}{2k-m-r}+k^n\phi_n.$$
Then the last equation can be written as
$$(m+r-k)\phi_{n}=k\phi_{n+1}.$$
After short calculations, we obtain
 $$\phi_n=\frac{\beta_1(2k-m-r)-k(k-m)}{k}\left(\frac{m+r-k}{k}\right)^{n-1}.$$
Hence,
$$\beta_n=\frac{(k-m)k^n}{2k-m-r}+\left(\beta_1(2k-m-r)-k(k-m)\right)\left(m+r-k\right)^{n-1}.$$
Thus
$$\beta_1=\frac{(k-m)k}{2k-m-r}+\left(\beta_1(2k-m-r)-k(k-m)\right) \ \ \Rightarrow \ \ \beta_1=\frac{(k-m)k}{2k-m-r}.$$
Then
$$\beta_n=\frac{(k-m)k^n}{2k-m-r}.$$

Note that $\alpha_n+\beta_n=k^n$, then
$$\alpha_n=\frac{k^n(k-r)}{2k-r-m}+\left(k(k-m)-\beta_1(2k-m-r)\right)\left(m+r-k\right)^{n-1}.$$
Since $$\beta_1=\frac{(k-m)k}{2k-m-r},$$
one gets
$$\alpha_n=\frac{k^n(k-r)}{2k-r-m}.$$
Consequently, it is easy to check that
 \begin{equation*}\label{new}\lim_{n\to\infty}\frac{(k-1)\alpha_n}{k^{n+1}-1}=\frac{(k-1)(k-r)}{k(2k-m-r)}\end{equation*}

  and
 \begin{equation}\label{new}\lim_{n\to\infty}\frac{(k-1)\beta_n}{k^{n+1}-1}=\frac{(k-1)(k-m)}{k(2k-m-r)}.\end{equation}
Then
\begin{equation}\label{new1}F_{ALT}(h)=-\frac{1}{\beta}\cdot \left[ \frac{(k-1)(k-r)\ln h-(k^2-(m+1)k+m)\ln l}{k(2k-m-r)}+\lim_{n\to \infty}\frac{(k-1)\ln \left(\sum_{i=0}^{|V_n|}\lambda^{C_{|V_n|}^i}\right)}{k^{n+1}-1}\right].\end{equation}
By AM-GM inequality
$$\sum_{i=0}^{|V_n|}\lambda^{C_{|V_n|}^i}\geq |V_n|\cdot \sqrt[|V_n|]{\lambda^{2^{|V_n|}}}=|V_n|\cdot\lambda^{2^{|V_n|}\cdot|V_n|^{-1}}.$$
Since $\ln x$ is an increasing function
 $$\ln \left(\sum_{i=0}^{|V_n|}\lambda^{C_{|V_n|}^i}\right)\geq \ln |V_n|+{2^{|V_n|}\cdot|V_n|^{-1}}\ln\lambda.$$
Then
$$\lim_{n\to \infty}\frac{(k-1)\ln \left(\sum_{i=0}^{|V_n|}\lambda^{C_{|V_n|}^i}\right)}{k^{n+1}-1}\geq \lim_{n\to \infty}\frac{\ln |V_n|}{|V_n|}+\lim_{n\to \infty}{2^{|V_n|}\cdot|V_n|^{-2}}\ln\lambda.$$
If $\lambda>1$ then
\begin{equation}\label{new2}\lim_{n\to \infty}\frac{(k-1)\ln \left(\sum_{i=0}^{|V_n|}\lambda^{C_{|V_n|}^i}\right)}{k^{n+1}-1}=\infty.\end{equation}
By (\ref{new}), (\ref{new1}) and (\ref{new2}) one gets
 \begin{equation}\label{law} F_{ALT}(h)=-\lim_{n\to \infty}\frac{1}{ \beta |V_n|}
\ln \left[h^{\alpha_n}l^{\beta_n}\left(\sum_{i=0}^{|V_n|}\lambda^{C_{|V_n|}^{i}}\right)\right]=-\infty.\end{equation}
Also, if $\lambda\in (0,1]$ then
$$0\leq \lim_{n\to \infty}\frac{(k-1)\ln \left(\sum_{i=0}^{|V_n|}\lambda^{C_{|V_n|}^i}\right)}{k^{n+1}-1}\leq \lim_{n\to \infty}\frac{\ln |V_n|}{|V_n|}=\lim_{n\to \infty}\ln \sqrt[|V_n|]{|V_n|}=0.$$

Namely, from (\ref{new1})
$$F_{ALT}(h)=-\frac{1}{\beta}\cdot \left[\frac{(k-1)(k-r)\ln h-(k^2-(m+1)k+m)\ln l}{k(2k-m-r)}\right].$$

Hence from above results and by Theorem 5 and Theorem 6 we can conclude the following theorem.

\textbf{Theorem 7.} a) Let $k\geq2$, $m+r\leq k-2$ ($m=r$) and $\lambda^{(1)}_{cr}=\left(\frac{k-2m}{k-2m-1}\right)^k\cdot\frac{1}{k-2m-1}$. Then the following statements are true
 \begin{itemize} \item if $\lambda_{cr}^{(1)}\in (0,1]$ and $\lambda\in [\lambda_{cr}^{(1)}, 1]$ (resp. $\lambda\in (1, +\infty)$) then free energies $F_{ALT}$ of b.c (\ref{e3}) is equal to
$$-\frac{1}{\beta}\cdot \left[\frac{(k-1)(k-r)\ln h-(k^2-(m+1)k+m)\ln l}{k(2k-m-r)}\right] (\textrm{resp.} \ -\infty).$$
\item if $\lambda_{cr}^{(1)}\in (1, \infty)$ then free energies $F_{ALT}$ equals $-\infty$.
\end{itemize}

b) Let $k\geq2$, $r+m=k-2$ and $C_k^{m+1}\lambda_{cr}\leq 2^k$. Then the following statements hold:
\begin{itemize}

\item  if $m\neq r$ and $\lambda=\lambda_{cr}$ then free energies $F_{ALT}$ equals $-\infty$. Also, if $\lambda>\lambda_{cr}$ then free energies $F_{ALT}$ equals $-\infty$.
\end{itemize}

\selectlanguage{english}

\end{document}